\documentclass{amsart}
\usepackage{amsthm}                 
\usepackage{amsmath}                
\usepackage{amsfonts,amssymb}       
\usepackage{graphicx}               
\usepackage{float}                  
\usepackage{rotating}               
\restylefloat{figure}               
\usepackage{mathrsfs}               
\usepackage{fancyhdr}               
\usepackage[breaklinks]{hyperref}   
\usepackage[numbers]{natbib}
\allowdisplaybreaks
\theoremstyle{plain}
\newtheorem{lem}{Lemma}
\newtheorem{prop}[lem]{Proposition}
\newtheorem{thm}[lem]{Theorem}

\theoremstyle{definition}

\newcommand{\R}{\mathbb R}

\newcommand{\Z}{\mathbb Z}
\newcommand{\N}{\mathbb N}

\newcommand{\Diff}{\mbox{\rm Diff}}

\newcommand{\id}{\text{\rm id}}

\newcommand{\Ad}{\text{\rm Ad}}
\newcommand{\ad}{\text{\rm ad}}

\renewcommand{\S}{\mathbb S}

\renewcommand{\phi}{\varphi}

\newcommand{\norm}[1]{\left|\!\left|#1\right|\!\right|}

\newcommand{\eps}{\varepsilon}

\newcommand{\ska}[2]{\left\langle #1,#2\right\rangle}

\newcommand{\bea}{\begin{eqnarray}}
\newcommand{\eea}{\end{eqnarray}}
\newcommand{\beq}{\begin{equation}}
\newcommand{\eeq}{\end{equation}}
\renewcommand{\phi}{\varphi}

\renewcommand{\autoref}[1]{\text{Eq.}~\eqref{#1}}
\begin{document}
\title{Curvature computations for a two-component Camassa-Holm equation with vorticity}
\author{Martin Kohlmann}
\address{Dr.\ Martin Kohlmann, Goerdelerstra{\ss}e 36, 38228 Salzgitter, Germany}
\email{martin\_kohlmann@web.de}
\begin{abstract} In the present paper, a two-component Camassa-Holm (2CH) system with vorticity is studied as a geodesic flow on a suitable Lie group. The paper aims at presenting various details of the geometric formalism and a major result is the computation of the sectional curvature $K$ of the underlying configuration manifold. As a further result, we show that there are directions for which $K$ is strictly positive and bounded away from zero.
\end{abstract}
\maketitle
\tableofcontents
\section{Introduction}\label{sec_intro}
As a mathematical model for two-dimensional shallow water waves with constant vorticity, the following generalized two-component Camassa-Holm (2CH) system has attracted a considerable amount of interest recently:
\beq\label{2CH}
\left\{
\begin{array}{rcl}
  m_t       & = & \alpha u_x -au_xm -um_x -\kappa\rho\rho_x,\\
  \rho_t    & = & -u\rho_x-(a-1)u_x\rho,\\
  \alpha_t  & = & 0.
\end{array}
\right.
\eeq
Here $a\in\R\backslash\{1\}$, $\alpha$ is a constant, $\kappa>0$ and $m=Au$ with $A$ denoting the Fourier multiplication operator $A=(1-\partial_x^2)^s$ for $s\geq 1$. The functions $u$ and $\rho$ depend on time $t$ and a spatial variable $x\in\S\simeq\mathbb R/\Z$. A derivation of the system \eqref{2CH} with $s=1$ by means of formal asymptotic methods applied to the full-governing equations for two-dimensional water waves with constant vorticity is the subject of the paper \cite{EHKL14}. A special case of the system \eqref{2CH} is a one-parameter family of evolution equations obtained for $\alpha=0$, $s=1$ and $\rho\equiv 0$ and it is also called the $b$-equation \cite{DGH03,DHH02,HS03} (here with the parameter $a\neq 1$). To further special cases of this family are the Camassa-Holm (CH) equation ($a=2$)
\beq\label{CH}
u_t-u_{txx}=2u_xu_{xx}-3uu_x+uu_{xxx},
\eeq
cf.\ \cite{CH93}, and the Degasperis-Procesi (DP) equation ($a=3$)
$$
u_t-u_{txx}=3u_xu_{xx}-4uu_x+uu_{xxx},
$$
cf.\ \cite{DP99}. Apart from the fact that both the CH equation and the DP equation are of hydrodynamical relevance, see e.g.\ \cite{DGH03,I07}, they share some interesting mathematical properties: the $b$-equation is integrable only if $b\in\{2,3\}$ in the sense that for $b\in\{2,3\}$ there exists a bi-Hamiltonian formulation and a Lax pair representation \cite{CH93,DHH02}. Moreover, both the CH and the DP equation allow for finite-time solutions that can be interpreted as breaking waves \cite{CE00} or shock waves \cite{ELY07} as well as for global solutions \cite{EY08} and peculiar traveling wave solutions \cite{L05CH,L05DP}.\ The 2CH equation without vorticity ($\alpha=0$, $a=2$ and $s=1$ in \eqref{2CH}) has been the subject of \cite{CI08,ELeY07} where the authors proved local-in-time well-posedness by applying Kato's semigroup approach, discussed blow-up and established its integrable structure. In addition, the hydrodynamical relevance of the 2CH system without vorticity as a two-component extension of the CH equation is presented.

A remarkable property of the $b$-family equations is that they allow for a geometric reformulation on the diffeomorphism group of the circle. The group of all smooth and orientation-preserving diffeomorphisms $\S\to\S$, denoted as $\Diff^\infty(\S)$, is a Fr\'echet Lie group that can be equipped with an affine connection $\nabla$ so that the $b$-equation is the geodesic equation on $\Diff^\infty(\S)$ with respect to the connection $\nabla$, see, e.g.\ \cite{EK10,K09}. Furthermore, the geodesic flow is the minimizer of a length functional if and only if $b=2$; in this case the resulting equation, precisely the CH equation, is a metric Euler equation with respect to the $H^1$ metric. The DP equation belongs to the class of non-metric Euler equations \cite{ES10,K09}. An analogous geometric framework has been established for the 2CH equation without vorticity in \cite{EKL11,HT09} where the authors showed that it can be recast as a geodesic flow on the semidirect product $\Diff^\infty(\S)\circledS C^\infty(\S)$ equipped with the $H^1$ metric for the first component plus the $L_2$ metric for the second component.

The geometric theory for evolution equations arising in hydrodynamics is not only a technical game: There are various applications of the geometric picture concerning the features of solutions to report on. In \cite{EK10}, the authors make use of the geometric reformulation to establish the well-posedness of the periodic $b$-equation on a scale of Sobolev spaces that are used as a Banach space approximation for $C^\infty(\S)$. In \cite{CK02}, finite time solutions to the CH equation are related to a breakdown of the geodesic flow. In particular, computations of the sectional curvature have been performed as the sign of the sectional curvature of the underlying configuration manifold has implications for the stability of its geodesics. The curvature tensor for the CH equation on the diffeomorphism group of $\S$ has been computed in \cite{M02} where the author also shows that the non-normalized sectional curvature $S_{\text{CH}}$ is positive on an infinite-dimensional subspace containing the vector $\cos kx$ and that the normalized sectional curvature is bounded away from zero on this subspace. Moreover,
$$S_{\text{CH}}(u,v)=\ska{\Gamma_{\text{CH}}(u,v)}{\Gamma_{\text{CH}}(u,v)}_{H^1}-\ska{\Gamma_{\text{CH}}(u,u)}{\Gamma_{\text{CH}}(v,v)}_{H^1}$$
with $\Gamma_{\text{CH}}$ denoting the Christoffel operator for the CH equation. Similarly, it could be shown in \cite{EKL11} that the non-normalized sectional curvature $S_{\text{2CH}}$ for the 2CH equation without vorticity is strictly positive in directions spanned by vectors of the type $(\cos kx,\cos lx)$, that the normalized sectional curvature is bounded away from zero in directions spanned by vectors of the type $(0,\cos kx)$ and that
$$S_{\text{2CH}}(u,v)=\ska{\Gamma_{\text{2CH}}(u,v)}{\Gamma_{\text{2CH}}(u,v)}_{H^1\oplus L_2}-\ska{\Gamma_{\text{2CH}}(u,u)}{\Gamma_{\text{2CH}}(v,v)}_{H^1\oplus L_2},$$
with $\Gamma_{\text{2CH}}$ denoting the Christoffel operator for the 2CH equation. We refer the reader to \cite{BKP15,KLMP13} for further curvature computations for related equations of hydrodynamical relevance.

In the present paper, we first focus in detail on the geometric picture for \autoref{2CH} with the fixed parameters $a=2$ and $\kappa,s=1$. The fact that \eqref{2CH} represents geodesic motion on the Lie group $(\Diff^{\infty}(\S)\circledS C^\infty(\S))\times\R$ with a suitable right-invariant metric $\ska{\cdot}{\cdot}_{\mathbb A}$ induced by an operator $\mathbb A$ has been established in \cite{EHKL14} where the authors compute the adjoint of the adjoint action on the Lie algebra with respect to $\ska{\cdot}{\cdot}_{\mathbb A}$ in order to identify \autoref{2CH} with the geodesic equation $U_t=-\ad_U^*U$, where $U=(u,\rho,\alpha)$. However, the various analogies and their consequences when comparing the geometric picture for \autoref{2CH} with the geometric picture for the rigid body motion pioneered by Arnold \cite{A66} in 1966 have not been work out to the best of the author's knowledge. Section~\ref{sec_geometry} of the present work has the goal to provide some further aspects of the geometric theory for the 2CH equation \eqref{2CH}. In Section~\ref{sec_curv}, we present the following main theorem on the sectional curvature associated with \autoref{2CH}. It clearly gets in line with the above mentioned results on the sectional curvature for the CH equation and the 2CH equation with $\alpha=0$.
\thm\label{thm_curvature} Let $R$ denote the curvature tensor associated with the 2CH equation \eqref{2CH} on $(\Diff^{\infty}(\S)\circledS C^\infty(\S))\times\R$ and denote by $S(u,v)=\ska{R(u,v)v}{v}_{\mathbb A}$ the non-normalized sectional curvature at the identity. Then
\beq\label{seccurv}S(u,v)=\ska{\Gamma(u,v)}{\Gamma(u,v)}_{\mathbb A}-\ska{\Gamma(u,u)}{\Gamma(v,v)}_{\mathbb A}.\eeq
Moreover, $S(u,v)>0$ for all vectors of the form
$$u=\begin{pmatrix}\cos k_1x\\\cos k_2x\\\alpha\end{pmatrix},\quad v=\begin{pmatrix}\cos l_1x\\\cos l_2x\\1\end{pmatrix},\quad k_1,k_2,l_1,l_2\in 2\pi\N,
\quad \alpha\geq 6\max\{k_1^2,l_1^2\},$$
and the normalized sectional curvature
$$K(u,v)=\frac{S(u,v)}{\ska{u}{u}_{\mathbb A}\ska{v}{v}_{\mathbb A}-\ska{u}{v}_{\mathbb A}^2}$$
is bounded away from zero for fixed values $k_1\neq k_2$, $l_1\neq l_2$ and as $\alpha\to\infty$.
\endthm\rm
\section{The geometric formalism}\label{sec_geometry}
In 1966, Arnold \cite{A66} showed that the motion of a rigid body rotating around its center of mass is in fact geodesic motion on the group $G=SO(3)$. The configuration of the body at time $t$ is given by a rotation matrix $R(t)$ which maps the position of a particle in body coordinates to its spatial position. The quantities $\omega=\dot{R}R^{-1}$ and $\Omega=R^{-1}\dot{R}$ are elements of $\mathfrak{so}(3)$, the Lie algebra $\mathfrak g$ of $SO(3)$, and correspond to the angular velocity in spatial or body coordinates respectively. A moment of inertia tensor $\mathbb I$ maps the body velocity to its momentum $\Pi=\mathbb I\Omega$; the spatial momentum is denoted by $\pi=R\Omega$. As $\mathfrak{so}(3)$ and $\mathfrak{so}(3)^*$ are canonically identified with $\R^3$, the Adjoint and Co-Adjoint actions $\text{Ad}_R$, $\text{Ad}^*_R$ are maps $\R^3\to\R^3$ and they link the velocity and the momentum in the spatial and the body frame of reference. Euler's equations for the rigid body motion particularly imply that the spatial momentum is a conserved quantity.

Ebin and Marsden \cite{EM70} proved in 1970 that Arnold's formalism can also be applied to the motion of an ideal fluid for which the configuration space is the group of all volume-preserving diffeomorphisms of the fluid domain. A major difference to the rigid body motion is that the Riemannian metric on the diffeomorphism group is right-invariant whereas the geodesic equation for the rigid body is induced by a left-invariant metric. Details of this geometric approach have been elaborated in detail for the Camassa-Holm equation \eqref{CH} in \cite{K99,L07} where the authors show that \autoref{CH} with periodic boundary conditions is equivalent to a geodesic equation on the group $\Diff^\infty(\S)$ of all smooth and orientation-preserving diffeomorphisms on $\S$. In \cite{EKL11,HT09} the authors showed that the 2CH equation without vorticity ($\alpha=0$) allows for a geometric reformulation on the semidirect product group $\Diff^\infty(\S)\circledS C^\infty(\S)$. We refer the reader to Appendix A of the paper \cite{EKL11} where the analogy of Arnold's approach to the geometric picture for CH and 2CH without vorticity is explained.

In this section, we present in detail the geometric picture for \autoref{2CH} with the fixed parameters $a=2$ and $\kappa,s=1$.
\subsection{The Lie group}
Consider the Fr\'echet Lie group
$$C^\infty G:=(\Diff^\infty(\S)\circledS C^\infty(\S))\times\R$$
where $\Diff^\infty(\S)$ denotes the group of smooth and orientation-preserving diffeomorphisms of $\S:=\S^1\simeq\R/\Z$ and $\circledS$ denotes a semidirect product. Writing $\circ$ for the composition of functions, the group product on $C^\infty G$ is given by
$$(\phi_1,f_1,s_1)*(\phi_2,f_2,s_2)=(\phi_1\circ\phi_2,f_2+f_1\circ\phi_2,s_1+s_2),$$
for $(\phi_1,f_1,s_1),(\phi_2,f_2,s_2)\in \Diff^\infty(\S)\times C^\infty(\S)\times\R$. The neutral element on $C^\infty G$ is $(\id,0,0)$ and one easily checks that $(\phi,f,s)\in C^\infty G$ has the inverse
$$(\phi,f,s)^{-1}=(\phi^{-1},-f\circ\phi^{-1},-s).$$
Let $R_g$ and $L_g$ denote right and left translation on $C^\infty G$ and write $I_gh=L_gR_{g^{-1}}h$ for the inner automorphism. We observe that
$$I_{(\phi_1,f_1,s_1)}(\phi_2,f_2,s_2)=(\phi_1\circ\phi_2\circ\phi_1^{-1},(f_2-f_1)\circ\phi_1^{-1}+f_1\circ(\phi_2\circ\phi_1^{-1}),s_2).$$
%
%
Writing $T_g$ for the tangent map at $g\in C^\infty G$, we have that
\begin{align}
\text{Ad}_{(\phi_1,f_1,s_1)}(u_2,\rho_2,\alpha_2)&=
\left[T_{(\id,0,0)}I_{(\phi_1,f_1,s_1)}\right](u_2,\rho_2,\alpha_2)\nonumber\\
&=((u_2\phi_{1x})\circ\phi_1^{-1},(\rho_2+f_{1x}u_2)\circ\phi_1^{-1},\alpha_2)\nonumber
\end{align}
%
and
\begin{align}\text{ad}_{(u_1,\rho_1,\alpha_1)}(u_2,\rho_2,\alpha_2)&=
T_{(\id,0,0)}\left[\text{Ad}_{(\cdot)}(u_2,\rho_2,\alpha_2)\right](u_1,\rho_1,\alpha_1)\nonumber\\
&=(u_{1x}u_2-u_{2x}u_1,\rho_{1x}u_2-\rho_{2x}u_1,0)\nonumber\\
&=[(u_1,\rho_1,\alpha_1),(u_2,\rho_2,\alpha_2)];\nonumber
\end{align}
here, $[\cdot,\cdot]$ denotes the Lie bracket on the Lie algebra
$$C^\infty\mathfrak g := T_{(\id,0,0)}C^\infty G\simeq C^\infty(\S)\times C^\infty(\S)\times\R.$$
For the following considerations, it will be also important to note the trivialization
$$TC^\infty G\simeq(\Diff^\infty(\S)\times C^\infty(\S)\times\R)\times(C^\infty(\S)\times C^\infty(\S)\times\R).$$
For a smooth path $g(t)=(\varphi,f,s)(t)$ in $C^\infty G$, the associated Eulerian velocity $(u,\rho,\alpha)(t)=T_{g(t)}R_{g^{-1}(t)}g'(t)\in C^\infty\mathfrak g$ is given by
\beq\label{Eulerian_vel}(u,\rho,\alpha)(t)=(\phi'\circ\phi^{-1},f'\circ\phi^{-1},s')(t)\eeq
and
$$U_0(t):=T_{g(t)}L_{g^{-1}(t)}g'(t)=\left(\frac{\phi_t}{\phi_x},f_t-f_x\frac{\phi_t}{\phi_x},s'\right)$$
so that, for $U=(u_1,u_2,u_3)\in C^\infty\mathfrak g$,
$$\Ad_{(\phi,f,s)}U=\left(\Ad_\phi u_1,(u_2+f_xu_1)\circ\phi^{-1},u_3\right),$$
where $\Ad_\phi u_1=(u_1\phi_x)\circ\phi^{-1}$ is the Adjoint action with respect to $\Diff^\infty(\S)$.
\subsection{The right-invariant metric}
We define an inner product on $C^\infty\mathfrak g$ by setting
\beq\label{metric}
\ska{U}{V}_{(\id,0,0)} := \int_\S u_1v_1\,dx+\int_\S u_{1x}v_{1x}\,dx+\int_\S u_2v_2\,dx
-\frac{1}{2}\int_\S (u_1v_3+u_3v_1)\,dx+\frac{1}{2}u_3v_3,
\eeq
where $U=(u_1,u_2,u_3),V=(v_1,v_2,v_3)\in C^\infty\mathfrak g$. It is shown in \cite{EL15} that $\ska{\cdot}{\cdot}_{(\id,0,0)}$ is indeed positive definite. With the inertia operator $\mathbb A\colon C^\infty\mathfrak{g}\to (C^\infty\mathfrak{g})^*$ given by
$$\mathbb AU:=\left(Au_1-\frac{1}{2}u_3,u_2,\frac{1}{2}\left(u_3-\int_\S u_1\,dx\right)\right),$$
where $A=1-\partial_x^2$, we observe that
$$\ska{U}{V}_{\mathbb A}:=\int_\S(\mathbb AU)\cdot Vdx=\ska{U}{V}_{(\id,0,0)}$$
and that the associated quadratic form is equivalent to the Hilbert norm $\norm{u_1}_{H^1}^2+\norm{u_2}_{L_2}^2+|u_3|^2$. We define a right-invariant metric on $C^\infty G$ by setting
$$\ska{U}{V}_{(\phi,f,s)}=\ska{TR_{(\phi,f,s)^{-1}}U}{TR_{(\phi,f,s)^{-1}}V}_{\mathbb A}$$
for all $U,V\in T_{(\phi,f,s)}C^\infty G\simeq C^\infty(\S)\times C^\infty(\S)\times\R$. It is well-known that the right-invariant metric for the 2CH equation without vorticity depends smoothly on $(\phi,f)$, cf.\ \cite{EKL11}, and by the definition of $\ska{\cdot}{\cdot}_{(\phi,f,s)}$ and the fact that $C^\infty G$ is a Lie group, it is immediately clear that $\ska{\cdot}{\cdot}_{(\phi,f,s)}$ depends smoothly on ${(\phi,f,s)}$ so that $(C^\infty G,\ska{\cdot}{\cdot}_{\mathbb A})$ is indeed a (weak) Riemannian manifold.

The operator $\mathbb A$ maps the Eulerian velocity $(u,\rho,\alpha)$ to the momentum
$$\mu:=\mathbb A(u,\rho,\alpha)=\left(Au-\frac{1}{2}\alpha,\rho,\frac{1}{2}\left(\alpha-\int_\S u\,dx\right)\right).$$
Using the $L_2$-pairing to identify the regular part of $(C^\infty\mathfrak g)^*$ with $C^\infty(\S)\times C^\infty(\S)\times\R$ and that $\Ad^*_\phi m=(m\circ\phi)\phi_x^2$, $m=Au$, we observe that
\begin{align}
\ska{\mu}{\Ad_{(\phi,f,s)}V}&=\int_\S\left(Au-\frac{\alpha}{2}\right)\Ad_\phi v_1\,dx+\int_\S\rho\left[(f_xv_1+v_2)\circ\phi^{-1}\right]dx+\nonumber\\
&\quad+\int_S\left(\frac{\alpha}{2}-\frac{1}{2}\int_\S u\,dx\right)v_3\,dx\nonumber\\
&=\int_\S\left\{\left[(m\circ\phi)-\frac{\alpha}{2}\right]\phi_x^2+(\rho\circ\phi)f_x\phi_x\right\}v_1\,dx+\int_\S(\rho\circ\phi)\phi_xv_2\,dx\nonumber\\
&\quad+\int_\S\left(\frac{\alpha}{2}-\frac{1}{2}\int_\S u\,dx\right)v_3\,dx\nonumber
\end{align}
so that $\mu_0:=\Ad_{(\phi,f,s)}^*\mu$ is given by
$$\mu_0=
\left(\left[(m\circ\phi)-\frac{\alpha}{2}\right]\phi_x^2+(\rho\circ\phi)f_x\phi_x,(\rho\circ\phi)\phi_x,\frac{1}{2}\left(\alpha-\int_\S u\,dx\right)\right).$$
We also show that, in analogy to the rigid body motion, we now obtain a conservation law for the 2CH equation \eqref{2CH}.
\prop The quantity $\mu_0$ corresponding to the body momentum of the rigid body motion is a conserved quantity for the 2CH equation \eqref{2CH}, i.e.
$$\frac{d}{dt}\mu_0=0.$$
\endprop\rm
\proof A simple calculation shows that
\begin{align}
&\quad\frac{d}{dt}\left[\left(m\circ\phi-\frac{\alpha}{2}\right)\phi_x^2+(\rho\circ\phi)f_x\phi_x\right] \nonumber\\
&= [(m_t+2u_xm+um_x-\alpha u_x+\rho\rho_x)\circ\phi]\phi_x^2+[(\rho_t+u\rho_x+u_x\rho)\circ\phi] f_x\phi_x\nonumber\\
&=0.\nonumber
\end{align}
That the second component of $\mu_0$ is conserved follows from Lemma 6.1 in \cite{EHKL14} and the time derivative of the third component of $\mu_0$ is zero as
$$u_t=-\partial_x\left[A^{-1}(u^2+\tfrac{1}{2}u_x^2+\tfrac{1}{2}\rho^2-\alpha u)+\tfrac{1}{2}u^2\right],$$
cf.\ \autoref{2CHweak}.
\endproof
The adjoint of the operator $\text{ad}\colon C^\infty\mathfrak g\times C^\infty\mathfrak g\to C^\infty\mathfrak g$ has been computed in \cite{EHKL14,EL15} and writing
$\text{ad}_{(u_1,\rho_1,\alpha_1)}^*(u_2,\rho_2,\alpha_2)=(\tilde u,\tilde\rho,\tilde\alpha)$, one has
\begin{align}
\tilde u    &= A^{-1}(2u_{1x}Au_2+u_1Au_{2x}-\alpha_2u_{1x}+\rho_{1x}\rho_2)+\int_\S(u_{1x}Au_2+\rho_{1x}\rho_2)\,dx, \nonumber\\
\tilde\rho  &= (u_1\rho_2)_x, \nonumber\\
\tilde\alpha&= 2\int_\S(u_{1x}Au_2+\rho_{1x}\rho_2)\,dx.\nonumber
\end{align}
\subsection{The geodesic spray}
To obtain the weak formulation of \eqref{2CH}, we apply the operator $A^{-1}$ to the first equation and add the term $uu_x$ so that
\begin{align}
 u_t+uu_x   & = A^{-1}(\alpha u_x-2u_xm-um_x-\rho\rho_x+A(uu_x))\nonumber\\
            & = -A^{-1}\partial_x\left(u^2+\tfrac{1}{2}u_x^2+\tfrac{1}{2}\rho^2-\alpha u\right);\nonumber
\end{align}
observe that the terms including third order derivatives of $u$ cancel out on the right hand side. We may thus rewrite \eqref{2CH} as
\beq\label{2CHweak}
\left\{
\begin{array}{rcl}
  u_t + uu_x        & = & -A^{-1}\partial_x\left(u^2+\tfrac{1}{2}u_x^2+\tfrac{1}{2}\rho^2-\alpha u\right),\\
  \rho_t +u\rho_x   & = & -\rho u_x,\\
  \alpha_t          & = & 0.
\end{array}
\right.
\eeq
Writing $U=(u_1,u_2,u_3),V=(v_1,v_2,v_3)\in C^\infty\mathfrak g$ and $\tilde U=(u_1,u_2)$ and $\tilde V=(v_1,v_2)$, we introduce the bilinear operator
\begin{align}
\Gamma(U,V)
&:=\begin{pmatrix}-A^{-1}\partial_x\left(u_1v_1+\tfrac{1}{2}u_{1x}v_{1x}+\tfrac{1}{2}u_2v_2-\tfrac{1}{2}u_3v_1-\tfrac{1}{2}u_1v_3\right)\\
-\tfrac{1}{2}u_2v_{1x}-\tfrac{1}{2}u_{1x}v_2\\0\end{pmatrix}\nonumber\\
&=\begin{pmatrix}\Gamma^0(\tilde U,\tilde V)\\0\end{pmatrix}+\begin{pmatrix}\tfrac{1}{2}A^{-1}\partial_x(u_3v_1+u_1v_3)\\0\\0\end{pmatrix}\label{Gamma2CH}
\end{align}
where
\beq\label{Gamma2CHold}\Gamma^0(\tilde U,\tilde V):=\begin{pmatrix}-A^{-1}\partial_x\left(u_1v_1+\tfrac{1}{2}u_{1x}v_{1x}+\tfrac{1}{2}u_2v_2\right)\\
-\tfrac{1}{2}u_2v_{1x}-\tfrac{1}{2}u_{1x}v_2\end{pmatrix}\eeq
denotes the Christoffel operator for the 2CH equation without vorticity \cite{EKL11}. Again, the map $\Gamma\colon C^\infty\mathfrak g\times C^\infty\mathfrak g\to C^\infty\mathfrak g$ can be extended to a right-invariant bilinear operator $\Gamma_{(\phi,f,s)}\colon T_{(\phi,f,s)}C^\infty G\times T_{(\phi,f,s)}C^\infty G\to T_{(\phi,f,s)}C^\infty G$ by setting
$$\Gamma_{(\phi,f,s)}=TTR_{(\phi,f,s)}\circ\Gamma_{(\id,0,0)}\circ TR_{(\phi,f,s)^{-1}}.$$
The second order vector field $TC^\infty G\to TTC^\infty G$, $(g,U)\mapsto (g,U,U,\Gamma_g(U,U))$ is called the geodesic spray for the 2CH equation \eqref{2CH}.
\subsection{The geodesic equation}
Introducing Lagrangian variables $(\phi,f,s)(t)$ for \autoref{2CH} by setting
$$\phi'=u\circ\phi,\quad f'=\rho\circ\phi,\quad s'=\alpha$$
it follows that \eqref{2CHweak} is equivalent to the geodesic equation
\beq\label{geodesiceq}(\phi,f,s)''(t)=\Gamma_{(\phi,f,s)(t)}((\phi,f,s)'(t),(\phi,f,s)'(t)).\eeq
We finally give a rigorous proof of the fact that the geodesics $(\phi,f,s)(t)$ are in fact length-minimizing with respect to the metric $\ska{\cdot}{\cdot}_{\mathbb A}$.
\prop\label{thm_variation} Let $\gamma(t)\colon[0,T]\to C^\infty G$ denote the shortest path on $C^\infty G$ between fixed endpoints with respect to the metric $\ska{\cdot}{\cdot}_{\mathbb A}$. Then $(u,\rho,\alpha)(t)=T_{\gamma(t)}R_{\gamma^{-1}(t)}\gamma'(t)$ is a solution to the 2CH equation \eqref{2CH}, i.e.\ \eqref{2CH} is the Euler-Lagrange equation for the action functional
$$\mathfrak{a}(\gamma)=\frac{1}{2}\int_0^T\ska{\gamma'(t)}{\gamma'(t)}_{\gamma(t)}dt.$$
\endprop\rm
\proof Assume that $\gamma$ is a critical point in the space of paths for the functional $\mathfrak a$. Then
$$\left.\frac{d}{d\eps}\right|_{\eps=0}\mathfrak a(\gamma+\eps\eta)=0$$
for every path $\eta\colon[0,T]\to C^\infty G$ with endpoints at zero and such that $\gamma+\eps\eta$ is a small variation of $\gamma$ on $C^\infty G$. As
\begin{align}
\left.\frac{d}{d\eps}\right|_{\eps=0}\mathfrak a(\gamma+\eps\eta) & = \int_0^T\int_\S(\gamma_1'\circ\gamma_1^{-1})\left.\frac{d}{d\eps}\right|_{\eps=0}
\left[(\gamma_1'+\eps\eta_1')\circ(\gamma_1+\eps\eta_1)^{-1}\right]dt\,dx\nonumber\\
& + \int_0^T\int_\S(\gamma_1'\circ\gamma_1^{-1})_x\left.\frac{d}{d\eps}\right|_{\eps=0}
\left[(\gamma_1'+\eps\eta_1')_x\circ(\gamma_1+\eps\eta_1)^{-1}\right]_xdt\,dx\nonumber\\
& + \int_0^T\int_\S(\gamma_2'\circ\gamma_1^{-1})\left.\frac{d}{d\eps}\right|_{\eps=0}
\left[(\gamma_2'+\eps\eta_2')\circ(\gamma_1+\eps\eta_1)^{-1}\right]dt\,dx\nonumber\\
&-\frac{1}{2}\int_0^T\int_\S\left.\frac{d}{d\eps}\right|_{\eps=0}\left\{(\gamma_3'+\eps\eta_3')\cdot\left[(\gamma_1'+\eps\eta_1')
\circ(\gamma_1+\eps\eta_1)^{-1}\right]\right\}dt\,dx\nonumber\\
&+\frac{1}{4}\int_0^T\int_\S\left.\frac{d}{d\eps}\right|_{\eps=0}(\gamma_3'+\eps\eta_3')^2\,dt\,dx\nonumber\\
& = \mathcal I_1+\mathcal I_2+\mathcal I_3+\mathcal I_4+\mathcal I_5,\nonumber
\end{align}
where the prime indicates differentiation with respect to time, we can invoke a result presented in \cite{IK07} to conclude that
$$\mathcal I_1+\mathcal I_2=-\int_0^T\int_\S(\eta_1\circ\gamma_1^{-1})[u_t+3uu_x-u_{txx}-2u_xu_{xx}-uu_{xxx}]\,dt\,dx$$
with $u=\gamma_1'\circ\gamma_1^{-1}$. Differentiating the equation $\gamma_1\circ\gamma_1^{-1}=\id$ with respect to $t$ and $x$ yields expressions for the derivatives of $\gamma_1^{-1}$ that help us to conclude that
\begin{align}
&\left.\frac{d}{d\eps}\right|_{\eps=0}\left[(\gamma_2'+\eps\eta_2')\circ(\gamma_1+\eps\eta_1)^{-1}\right] =
\eta_2'\circ\gamma_1^{-1}-[(\partial_x\gamma_2')\circ\gamma_1^{-1}]\frac{\eta_1\circ\gamma_1^{-1}}{(\partial_x\gamma_1)\circ\gamma_1^{-1}}\nonumber\\
&=\partial_t(\eta_2\circ\gamma_1^{-1})+(\gamma_1'\circ\gamma_1^{-1})\partial_x(\eta_2\circ\gamma_1^{-1})
-(\eta_1\circ\gamma_1^{-1})\partial_x(\gamma_2'\circ\gamma_1^{-1}).\nonumber
\end{align}
Writing $\rho=\gamma_2'\circ\gamma_1^{-1}$, integration by parts and the boundary conditions for $\eta$ now show that
$$\mathcal I_3=-\int_0^T\int_\S(\eta_2\circ\gamma_1^{-1})[\rho_t+(\rho u)_x]\,dt\,dx -\int_0^T\int_\S(\eta_1\circ\gamma_1^{-1})\rho\rho_x\,dt\,dx.$$
Similar calculations yield that
\begin{align}
&\mathcal I_4 = -\frac{1}{2}\int_0^T\int_\S\big\{\eta_3'u+\gamma_3'\big[\partial_t(\eta_1\circ\gamma_1^{-1})
+(\gamma_1'\circ\gamma_1^{-1})\partial_x(\eta_1\circ\gamma_1^{-1})\nonumber\\
&-(\eta_1\circ\gamma_1^{-1})\partial_x(\gamma_1'\circ\gamma_1^{-1})\big]\big\}\,dt\,dx.\nonumber
\end{align}
As $\eta_3(t)\in\R$ for any $t\in[0,T]$, the first term in $\mathcal I_4$ vanishes due to the fact that
\begin{align}
\int_0^T\int_\S\eta_3'u\,dt\,dx & = -\int_0^T\int_\S\eta_3u_t\,dt\,dx\nonumber\\
& = -\int_0^T\eta_3\left(\int_\S u_t\,dx\right)dt\nonumber\\
& = \int_0^T\eta_3\left(\int_\S \partial_x\left[A^{-1}(u^2+\tfrac{1}{2}u_x^2+\tfrac{1}{2}\rho^2-\alpha u)+\tfrac{1}{2}u^2\right]dx\right)dt\nonumber\\
& = 0.\nonumber
\end{align}
With $\gamma_3'=\alpha$, the remaining term can be rewritten as
$$\mathcal I_4=\int_0^T\int_\S(\eta_1\circ\gamma_1^{-1})(\tfrac{1}{2}\alpha_t+\alpha u_x)\,dt\,dx.$$
We finally observe that
$$\mathcal I_5=\frac{1}{2}\int_0^T\int_\S\gamma_3'\eta_3'\,dt\,dx=-\frac{1}{2}\int_0^T\int_\S\alpha_t\eta_3\,dt\,dx.$$
Hence the critical point of the length functional $\mathfrak a$ is obtained from the equation
\begin{align}
&\int_0^T\int_\S(\eta_1\circ\gamma_1^{-1})[u_t+3uu_x-u_{txx}-2u_xu_{xx}-uu_{xxx}-\alpha u_x+\rho\rho_x-\tfrac{1}{2}\alpha_t]\,dt\,dx\nonumber\\
&+\int_0^T\int_\S(\eta_2\circ\gamma_1^{-1})[\rho_t+(\rho u)_x]\,dt\,dx+\frac{1}{2}\int_0^T\int_\S\alpha_t\eta_3\,dt\,dx=0.\nonumber
\end{align}
Since we can choose $\eta$ arbitrarily, we immediately obtain the system \eqref{2CH} from the above identity.
\endproof
\subsection{The affine connection}
The geodesic flow $(\phi,f,s)(t)$ is not only the minimizer of the length functional on $(C^\infty G,\ska{\cdot}{\cdot}_{\mathbb A})$, it is also the geodesic flow corresponding to the affine connection
\beq\label{connection}(\nabla_XY)(\phi,f,s):=DY(\phi,f,s)\cdot X(\phi,f,s)-\Gamma_{(\phi,f,s)}(X,Y)\eeq
where $X$ and $Y$ are smooth vector fields on $C^\infty G$. It follows immediately from the definition \eqref{connection} that $\nabla$ is a Riemannian covariant derivative as defined in \cite{L07}, i.e.,
\begin{itemize}
\item[(i)] $\nabla$ is $\mathbb R$-bilinear,
\item[(ii)] $X(\phi,f,s)=0$ implies that $(\nabla_XY)(\phi,f,s)=0$,
\item[(iii)] $\nabla_X(fY)=f\nabla_XY+X(f)Y$ for $f\in C^\infty(C^\infty G)$ and
\item[(iv)] $\nabla_XY-\nabla_YX=[X,Y]$
\end{itemize}
with the Lie bracket given locally by $$[X,Y](g)=DY(g)\cdot X(g)-DX(g)\cdot Y(g),\quad g\in C^\infty G.$$
\prop The Riemannian covariant derivative $(X,Y)\overset{\nabla}{\to}\nabla_XY$ defined in \eqref{connection} is compatible with the metric \eqref{metric}.
\endprop\rm
\proof Let $X(\phi,f,s)$, $Y(\phi,f,s)$ and $Z(\phi,f,s)$ be smooth vector fields on $C^\infty G$ and let $u=X(\phi,f,s)\circ\phi^{-1},v=Y(\phi,f,s)\circ\phi^{-1}$ and $w=Z(\phi,f,s)$. Let $\gamma(\eps)\subset C^\infty G$ be a smooth path such that $\gamma(0)=(\phi,f,s)$ and $\gamma'(0)=X(\phi,f,s)$. Then
\begin{align}
(X\ska{Y}{Z}_{\mathbb A})(\phi,f,s)&=\left.\frac{d}{d\eps}\right|_{\eps=0}\ska{\begin{pmatrix}Y_1(\gamma(\eps))\circ\gamma_1(\eps)^{-1}\\
Y_2(\gamma(\eps))\circ\gamma_1(\eps)^{-1}\\Y_3(\gamma(\eps))\end{pmatrix}}{\begin{pmatrix}Z_1(\gamma(\eps))\circ\gamma_1(\eps)^{-1}\\
Z_2(\gamma(\eps))\circ\gamma_1(\eps)^{-1}\\Z_3(\gamma(\eps))\end{pmatrix}}_{\mathbb A}\nonumber\\
&\hspace{-3cm}=\left.\frac{d}{d\eps}\right|_{\eps=0}\bigg\{\ska{\tilde Y(\gamma(\eps))\circ\gamma_1(\eps)^{-1}}{\tilde Z(\gamma(\eps))\circ
\gamma_1(\eps)^{-1}}_{H^1\oplus L_2}\nonumber\\
&\hspace{-2cm}-\frac{1}{2}\int_\S [Y_1(\gamma(\eps))\circ\gamma_1(\eps)^{-1}Z_3(\gamma(\eps))+
Z_1(\gamma(\eps))\circ\gamma_1(\eps)^{-1}Y_3(\gamma(\eps))]\,dx\nonumber\\
&\hspace{-2cm}+\frac{1}{2}Y_3(\gamma(\eps))Z_3(\gamma(\eps))\bigg\}\nonumber\\
&\hspace{-3cm}=\left.\frac{d}{d\eps}\right|_{\eps=0}\ska{\tilde Y(\gamma(\eps))\circ\gamma_1(\eps)^{-1}}{\tilde Z(\gamma(\eps))\circ
\gamma_1(\eps)^{-1}}_{H^1\oplus L_2}\nonumber\\
&\hspace{-2cm}-\frac{1}{2}\int_\S[w_3(DY_1\cdot X)\circ\phi^{-1}+v_1DZ_3\cdot X-u_1v_{1x}w_3]\,dx\nonumber\\
&\hspace{-2cm}-\frac{1}{2}\int_\S[v_3(DZ_1\cdot X)\circ\phi^{-1}+w_1DY_3\cdot X-u_1v_3w_{1x}]\,dx\nonumber\\
&\hspace{-2cm}+\frac{1}{2}(w_3 DY_3\cdot X+v_3 DZ_3\cdot X).\nonumber
\end{align}
On the other hand, we have
\begin{align}
\ska{\nabla_XY}{Z}_{(\phi,f,s)}&=\ska{\begin{pmatrix}(D\tilde Y\cdot X)\circ\phi^{-1}-\Gamma^0(\tilde u,\tilde v)\\DY_3\cdot X\end{pmatrix}-\frac{1}{2}\begin{pmatrix}A^{-1}\partial_x(u_3v_1+u_1v_3)\\0\\0\end{pmatrix}}{w}_{\mathbb A}\nonumber\\
&\hspace{-2.3cm}=\ska{(D\tilde Y\cdot X)\circ\phi^{-1}-\Gamma^0(\tilde u,\tilde v)}{\tilde w}_{H^1\oplus L_2}\nonumber\\
&\hspace{-2.3cm}-\frac{1}{2}\int_\S
[w_1DY_3\cdot X+w_3(DY_1\cdot X)\circ\phi^{-1}-w_{1x}(u_3v_1+u_1v_3)]\,dx+\frac{1}{2}w_3DY_3\cdot X\nonumber
\end{align}
where we have used that integrals of the type
$$\int_\S\Gamma^0(\tilde u,\tilde v)w_3\,dx=\int_\S A^{-1}(u_1v_1+\tfrac{1}{2}u_{1x}v_{1x}+\tfrac{1}{2}u_2v_2)w_{3x}\,dx=0$$
vanish. Clearly,
\begin{align}
\ska{\nabla_XZ}{Y}_{(\phi,f,s)}&=\ska{(D\tilde Z\cdot X)\circ\phi^{-1}-\Gamma^0(\tilde u,\tilde w)}{\tilde v}_{H^1\oplus L_2}\nonumber\\
&\hspace{-2.3cm}-\frac{1}{2}\int_\S
[v_1DZ_3\cdot X+v_3(DZ_1\cdot X)\circ\phi^{-1}-v_{1x}(u_3w_1+u_1w_3)]\,dx+\frac{1}{2}v_3DZ_3\cdot X.\nonumber
\end{align}
Applying \cite[Prop.\ 3.1]{EKL11} completes the proof of the proposition.
\endproof
\subsection{Summary and conclusions}
In the following tabular, we summarize some unifying features of the approach pioneered by V.I.\ Arnold by comparing the geometric quantities for the rigid body motion with the corresponding quantities that have been presented in this section for the 2CH equation with vorticity.
\begin{center}
\begin{scriptsize}
\begin{tabular}{|l|c|c|}
  \hline
   & Rigid body & 2CH with vorticity\\
  \hline
  configuration space   & $SO(3)$               &$C^\infty G=(\Diff^\infty(\S)\circledS C^\infty(\S))\times\R$\\
  Lie algebra           & $\mathfrak{so}(3)$    & $C^\infty(\S)\times C^\infty(\S)\times\R$\\
  material velocity     & $\dot R(t)$           & $(\phi,f,s)'(t)$\\
  spatial velocity      & $\omega=\dot RR^{-1}$ & $(u,\rho,\alpha) = (\phi'\circ\phi^{-1},f'\circ\varphi^{-1},s') $ \\
  body velocity         & $\Omega=R^{-1}\dot R$ & $U_0=(\tfrac{\phi_t}{\phi_x},f_t-f_x\tfrac{\phi_t}{\phi_x},s')$ \\
  inertia operator      & $\mathbb I$           & $\mathbb A$ \\
  spatial momentum      & $\pi=R\Pi$            & $ \mu=(Au-\tfrac{\alpha}{2},\rho,\tfrac{\alpha}{2}-\tfrac{1}{2}\int_\S u\,dx)$ \\
  body momentum         & $\Pi=\mathbb I\Omega$ & $ \mu_0=\begin{pmatrix}[(m\circ\phi)-\tfrac{\alpha}{2}]\phi_x^2+(\rho\circ\phi)f_x\phi_x\\(\rho\circ\phi)\phi_x\\\tfrac{\alpha}{2}-\tfrac{1}{2}\int_\S u\,dx\end{pmatrix}$  \\
  spatial velocity (Ad) & $\omega = \Ad_{R}\Omega$ & $(u,\rho,\alpha) = \Ad_{(\varphi,f,s)}U_0$\\
  body momentum (Ad*)   & $\Pi=\Ad^*_R\pi$      & $\mu_0 =\Ad_{(\varphi,f,s)}^*\mu$ \\
  momentum conservation & $\pi = \text{const.}$ & $\mu_0 = \text{const.}$ \\
  Lie bracket (ad)      & $[A,B]=AB-BA$         & $[(u_1,\rho_1,\alpha_1),(u_2,\rho_2,\alpha_2)]                                         =\begin{pmatrix}u_{1x}u_2-u_{2x}u_1\\\rho_{1x}u_2-\rho_{2x}u_1\\0\end{pmatrix}$\\
  $\text{ad}^*$         & $\text{ad}_A^*B=[B,A]$& $\text{ad}_{(u_1,\rho_1,\alpha_1)}^*(u_2,\rho_2,\alpha_2)$\\
  &&                                              $=
  \begin{pmatrix}A^{-1}(2u_{1x}Au_2+u_1Au_{2x}-\alpha_2u_{1x}+\rho_{1x}\rho_2)\\+\int_\S(u_{1x}Au_2+\rho_{1x}\rho_2)\,dx\\(u_1\rho_2)_x\\
  2\int_\S(u_{1x}Au_2+\rho_{1x}\rho_2)\,dx\end{pmatrix}$\\
  \hline
\end{tabular}
\end{scriptsize}
\end{center}
The geometric theory is not only aesthetically appealing but also helps to understand some important features of the solutions to the 2CH equation:

The authors of \cite{EHKL14} showed that the geodesic spray $((\phi,f,s),U,U,\Gamma_{(\phi,f,s)}(U,U))$ is smooth as a map $TH^sG\to TTH^sG$, for $s>5/2$, where $H^sG$ denotes the group $(\Diff^s(\S)\circledS H^{s-1}(\S))\times\R$ and $\Diff^s(\S)$ is the group of all orientation-preserving $H^s$ diffeomorphisms $\S\to\S$. The groups $H^sG$ are only topological groups (but not Lie groups), instead they are Banach manifolds (and not Fr\'echet manifolds) so that the Picard-Lindel\"of Theorem can be applied to conclude the existence of a local-in-time solution $(\phi,f,s)(t)$ to the geodesic equation \eqref{geodesiceq} for any pair of initial values $(u_0,\rho_0)\in H^s(\S)\times H^{s-1}(\S)$ and $\alpha\in\R$. A Hilbert approximation of $C^\infty G$ by the groups $H^sG$ then shows that \autoref{geodesiceq} also possesses a unique non-extendable solution with smooth dependence on the initial data in the smooth category. As $C^\infty G$ is a Lie group, composition and inversion are smooth maps so that the relation \eqref{Eulerian_vel} immediately implies that \autoref{2CH} possesses a unique maximal solution $(u,\rho)(t)\in C^\infty(\S)\times C^\infty(\S)$, $t\in J$, for any initial datum $(u_0,\rho_0)\in C^\infty(\S)\times C^\infty(\S)$ and any $\alpha\in\R$.

By Theorem~6.5 of \cite{EHKL14}, the solution $(u,\rho)(t)$ exists for all $t\geq 0$ provided $\norm{u_x(t)}_\infty$ is bounded on any bounded subinterval of $J$.
\section{The sectional curvature}\label{sec_curv}
In this section, we present some curvature computations providing a proof of Theorem~\ref{thm_curvature}.

We begin with the term $\ska{\Gamma(u,v)}{\Gamma(u,v)}_{\mathbb A}$ on the right hand side of \eqref{seccurv} which we intend to rewrite as $\ska{\Gamma^0(\tilde u,\tilde v)}{\Gamma^0(\tilde u,\tilde v)}_{\tilde{\mathbb A}}$ plus additional terms; again $\Gamma^0$ denotes the spray for the 2CH equation without vorticity, cf.\ \autoref{Gamma2CHold}, $\tilde u=(u_1,u_2)$, $\tilde v=(v_1,v_2)$ and $\tilde{\mathbb A}=\text{diag}(A,1)$. By \eqref{Gamma2CH} and the definition of the metric \eqref{metric}, we have that
\begin{align}
\ska{\Gamma(u,v)}{\Gamma(u,v)}_{\mathbb A} &= \ska{\Gamma^0(\tilde u,\tilde v)}{\Gamma^0(\tilde u,\tilde v)}_{\tilde{\mathbb A}} +
\ska{\Gamma^0(\tilde u,\tilde v)}{\begin{pmatrix}A^{-1}\partial_x(u_3v_1+u_1v_3)\\0\end{pmatrix}}_{\tilde{\mathbb A}}\nonumber\\
&\quad+\frac{1}{4}\ska{A^{-1}\partial_x(u_3v_1+u_1v_3)}{A^{-1}\partial_x(u_3v_1+u_1v_3)}_A\nonumber\\
&=\ska{\Gamma^0(\tilde u,\tilde v)}{\Gamma^0(\tilde u,\tilde v)}_{\tilde{\mathbb A}}\nonumber\\
&\quad-\int_\S\partial_x(u_1v_1+\tfrac{1}{2}u_{1x}v_{1x}+\tfrac{1}{2}u_2v_2)A^{-1}\partial_x(u_3v_1+u_1v_3)\,dx\nonumber\\
&\quad+\frac{1}{4}\int_\S\partial_x(u_3v_1+u_1v_3)A^{-1}\partial_x(u_3v_1+u_1v_3)\,dx.\label{RHSterm1}
\end{align}
The second term on the right hand side of \eqref{seccurv} is computed similarly and we find that
\begin{align}
\ska{\Gamma(u,u)}{\Gamma(v,v)}_{\mathbb A}
&=\ska{\Gamma^0(\tilde u,\tilde u)}{\Gamma^0(\tilde v,\tilde v)}_{\tilde{\mathbb A}}\nonumber\\
&\quad-\int_\S\partial_x(u_1^2+\tfrac{1}{2}u_{1x}^2+\tfrac{1}{2}u_2^2)A^{-1}\partial_x(v_1v_3)\,dx\nonumber\\
&\quad-\int_\S\partial_x(v_1^2+\tfrac{1}{2}v_{1x}^2+\tfrac{1}{2}v_2^2)A^{-1}\partial_x(u_1u_3)\,dx\nonumber\\
&\quad+\int_\S\partial_x(u_1u_3)A^{-1}\partial_x(v_1v_3)\,dx.\label{RHSterm2}
\end{align}
Let $u,v,w\in T_{p}C^\infty G$ be three tangent vectors at a point $p \in C^\infty G$.
The curvature tensor $R$ for $(C^\infty G,\ska{\cdot}{\cdot}_{\mathbb A})$ is given locally by
\beq\label{R}
R_p(u, v)w = D_1\Gamma_{p}(w, u)v - D_1\Gamma_{p}(w, v)u + \Gamma_{p}(\Gamma_{p}(w, v), u) - \Gamma_{p}(\Gamma_{p}(w, u), v),
\eeq
cf.\ \cite{Lang2}, where $D_1$ denotes differentiation with respect to
$p$:
$$D_1\Gamma_{p}(w, u)v = \frac{d}{d\epsilon}\bigg|_{\epsilon = 0} \Gamma_{p + \epsilon v}(w, u).$$
We apply \autoref{R} at $p=\id$ and with $w=v$ in order to rewrite $S(u,v)=\ska{R(u,v)v}{u}_{\mathbb A}$ as a sum of terms involving $\Gamma^0$ and first and second components of $u$ and $v$ plus additional terms involving the third components. Therefore, we make use of the identities
$$\left.\frac{d}{d\eps}\right|_{\eps=0} u_1\circ(\id+\eps v_1)^{-1} =-u_{1x}v_1$$
and
$$D_1\Gamma^0(\tilde w,\tilde u)\tilde v=-\Gamma^0(\tilde w_xv_1,\tilde u)-\Gamma^0(\tilde u_xv_1,\tilde w)+\Gamma^0(\tilde w,\tilde u)_xv_1,$$
see the proof of \cite[Prop.~5.1]{EKL11}, to infer that
\begin{align}
D_1\Gamma(v,u)v &= \begin{pmatrix}-\Gamma^0(\tilde v_xv_1,\tilde u)-\Gamma^0(\tilde u_xv_1,\tilde v)+\Gamma^0(\tilde v,\tilde u)_xv_1\\0\end{pmatrix}\nonumber\\
&\quad+\begin{pmatrix}-\tfrac{1}{2}A^{-1}\partial_x(u_{1x}v_1v_3+u_3v_1v_{1x})\\0\\0\end{pmatrix}\nonumber\\
&\quad+\begin{pmatrix}\tfrac{1}{2}v_1A^{-1}\partial_x^2(u_1v_3+u_3v_1)\\0\\0\end{pmatrix}\nonumber
\end{align}
and similarly that
\begin{align}
D_1\Gamma(v,v)u &= \begin{pmatrix}-2\Gamma^0(\tilde v_xu_1,\tilde v)+\Gamma^0(\tilde v,\tilde v)_xu_1\\0\end{pmatrix}\nonumber\\
&\quad+\begin{pmatrix}-A^{-1}\partial_x(u_1v_{1x}v_3)+u_1A^{-1}\partial_x^2(v_1v_3)\\0\\0\end{pmatrix}.\nonumber
\end{align}
Using once more the definition \eqref{Gamma2CH}, we also find that
\begin{align}
\Gamma(\Gamma(v,v),u) &= \begin{pmatrix}\Gamma^0(\Gamma^0(\tilde v,\tilde v),\tilde u)\\0\end{pmatrix}+\begin{pmatrix}\tfrac{1}{2}A^{-1}\partial_x(u_3\Gamma^0(\tilde v,\tilde v)_1)\\0\\0\end{pmatrix}\nonumber\\
&\hspace{-.5cm}+\begin{pmatrix}-A^{-1}\partial_x(u_1A^{-1}\partial_x(v_1v_3)+\frac{1}{2}u_{1x}A^{-1}\partial_x^2(v_1v_3)-\frac{1}{2}u_3A^{-1}\partial_x(v_1v_3))
\\-\tfrac{1}{2}u_2A^{-1}\partial_x^2(v_1v_3)\\0\end{pmatrix}\nonumber
\end{align}
and similarly that
\begin{align}
\Gamma(\Gamma(v,u),v) &= \begin{pmatrix}\Gamma^0(\Gamma^0(\tilde v,\tilde u),\tilde v)\\0\end{pmatrix}+\begin{pmatrix}\tfrac{1}{2}A^{-1}\partial_x(v_3\Gamma^0(\tilde v,\tilde u)_1)\\0\\0\end{pmatrix}\nonumber\\
&\hspace{-.05cm}+\begin{pmatrix}-\frac{1}{2}A^{-1}\partial_x(v_1A^{-1}\partial_x(u_1v_3+u_3v_1)+\frac{1}{2}v_{1x}A^{-1}\partial_x^2(u_1v_3+u_3v_1)\\
-\frac{1}{2}v_3A^{-1}\partial_x(u_1v_3+u_3v_1))\\[-.25cm]
\\-\tfrac{1}{4}v_2A^{-1}\partial_x^2(u_1v_3+u_3v_1)\\0\end{pmatrix}.\nonumber
\end{align}
As shown in the proof of \cite[Prop.~5.1]{EKL11}
\begin{align}\nonumber
&\ska{\Gamma^0(\tilde u,\tilde v)}{\Gamma^0(\tilde u,\tilde v)}_{\tilde{\mathbb A}}-\ska{\Gamma^0(\tilde u,\tilde u)}{\Gamma^0(\tilde v,\tilde v)}_{\tilde{\mathbb A}} =
\ska{\Gamma^0(\tilde v,\tilde u)_xv_1}{\tilde u}_{\tilde{\mathbb A}} - \ska{\Gamma^0(\tilde v,\tilde v)_xu_1}{\tilde u}_{\tilde{\mathbb A}}\nonumber\\
&\hspace{3cm}-\ska{\Gamma^0(\tilde v_xv_1,\tilde u)}{\tilde u}_{\tilde{\mathbb A}}-\ska{\Gamma^0(\tilde v,\tilde u_xv_1)}{\tilde u}_{\tilde{\mathbb A}}+2\ska{\Gamma^0(\tilde v_xu_1,\tilde v)}{\tilde u}_{\tilde{\mathbb A}}\nonumber\\
&\hspace{3cm}+\ska{\Gamma^0(\Gamma^0(\tilde v,\tilde v),\tilde u)}{\tilde u}_{\tilde{\mathbb A}}-\ska{\Gamma^0(\Gamma^0(\tilde v,\tilde u),\tilde v)}{\tilde u}_{\tilde{\mathbb A}}\nonumber
\end{align}
so that, using once more the definition \eqref{metric}, we obtain that
\beq
S(u,v)=\ska{\Gamma^0(\tilde u,\tilde v)}{\Gamma^0(\tilde u,\tilde v)}_{\tilde{\mathbb A}}-\ska{\Gamma^0(\tilde u,\tilde u)}{\Gamma^0(\tilde v,\tilde v)}_{\tilde{\mathbb A}} +\mathcal J_1+\mathcal J_2+\mathcal J_3\label{RHSSuv}
\eeq
where
\begin{align}
\mathcal J_1 &= -\frac{1}{2}\int_\S (u_1v_3+u_3v_1)\Gamma^0(\tilde v,\tilde u)_{1x}\,dx+\frac{1}{2}\int_\S(Au_1)v_1A^{-1}\partial_x^2(u_1v_3+u_3v_1)\,dx\nonumber\\
&-\frac{1}{4}\int_S u_3v_1A^{-1}\partial_x^2(u_1v_3+u_3v_1)\,dx+\frac{1}{4}\int_S u_2v_2A^{-1}\partial_x^2(u_1v_3+u_3v_1)\,dx\nonumber\\
&-\frac{1}{2}\int_\S u_{1x}\big[v_1A^{-1}\partial_x(u_1v_3+u_3v_1)+\tfrac{1}{2}v_{1x}A^{-1}\partial_x^2(u_1v_3+u_3v_1)\nonumber\\
&\hspace{2cm}-\tfrac{1}{2}v_3A^{-1}\partial_x(u_1v_3+u_3v_1)\big]\,dx\nonumber\\
&+\frac{1}{2}\int_\S u_{1x}(u_{1x}v_1v_3+v_1v_{1x}u_3)\,dx,\nonumber\\
\mathcal J_2 &=\int_\S u_1u_3\Gamma^0(\tilde v,\tilde v)_{1x}\,dx+\frac{1}{2}\int_\S u_1u_3A^{-1}\partial_x^2(v_1v_3)\,dx-\frac{1}{2}\int_\S u_{1x}u_3A^{-1}\partial_x(v_1v_3)\,dx,\nonumber\\
\mathcal J_3 &=\int_\S u_{1x}\left[u_1A^{-1}\partial_x(v_1v_3)+\tfrac{1}{2}u_{1x}A^{-1}\partial_x^2(v_1v_3)\right]dx-\frac{1}{2}\int_\S u_2^2A^{-1}\partial_x^2(v_1v_3)\,dx\nonumber\\
&-\int_\S(Au_1)u_1A^{-1}\partial_x^2(v_1v_3)\,dx-\int_\S u_1u_{1x}v_{1x}v_3\,dx.\nonumber
\end{align}
Now the proof of formula \eqref{seccurv} is completed by verifying that all the terms including third components on the right hand sides of \eqref{RHSterm1} and \eqref{RHSterm2} are equal to the third component terms $\mathcal J_1+\mathcal J_2+\mathcal J_3$ on the right hand side of \eqref{RHSSuv}. Using that $u_{3x}=v_{3x}=0$, that $\partial_x^2A^{-1}=A^{-1}\partial_x^2=-1+A^{-1}$ and integration by parts, some tedious computations which are omitted for the convenience of the reader show that indeed
\begin{align}
\mathcal J_1+\mathcal J_2+\mathcal J_3 &= -\int_\S\partial_x(u_1v_1+\tfrac{1}{2}u_{1x}v_{1x}+\tfrac{1}{2}u_2v_2)A^{-1}\partial_x(u_3v_1+u_1v_3)\,dx\nonumber\\
&\quad+\frac{1}{4}\int_\S\partial_x(u_3v_1+u_1v_3)A^{-1}\partial_x(u_3v_1+u_1v_3)\,dx\nonumber\\
&\quad+\int_\S\partial_x(u_1^2+\tfrac{1}{2}u_{1x}^2+\tfrac{1}{2}u_2^2)A^{-1}\partial_x(v_1v_3)\,dx\nonumber\\
&\quad+\int_\S\partial_x(v_1^2+\tfrac{1}{2}v_{1x}^2+\tfrac{1}{2}v_2^2)A^{-1}\partial_x(u_1u_3)\,dx\nonumber\\
&\quad-\int_\S\partial_x(u_1u_3)A^{-1}\partial_x(v_1v_3)\,dx.\nonumber
\end{align}
We now let
$$u=\begin{pmatrix}\cos k_1x\\\cos k_2x\\\alpha\end{pmatrix},\quad v=\begin{pmatrix}\cos l_1x\\\cos l_2x\\\beta\end{pmatrix},$$
for $k_1,k_2,l_1,l_2\in 2\pi\N$, and apply the identities
\begin{align*}
\cos\xi_1\cos\xi_2&=\frac{1}{2}(\cos(\xi_1-\xi_2)+\cos(\xi_1+\xi_2)),
    \\
\sin\xi_1\sin\xi_2&=\frac{1}{2}(\cos(\xi_1-\xi_2)-\cos(\xi_1+\xi_2)),
    \\
\sin\xi_1\cos\xi_2&=\frac{1}{2}(\sin(\xi_1-\xi_2)+\sin(\xi_1+\xi_2)),
\end{align*}
and
\begin{align*}
 A^{-1}\cos\xi x= &\;\frac{1}{1+\xi^2}\cos \xi x, \qquad \xi \in \R,
    \\
 \int_0^1\cos(\xi_1 x)\cos(\xi_2 x) dx= &\;\frac{1}{2}\left(\delta_{\xi_1,\xi_2}+\delta_{\xi_1,-\xi_2}\right), \qquad \xi_1,\xi_2\in 2\pi\Z,
    \\
 \int_0^1\sin(\xi_1 x)\sin(\xi_2 x)dx = &\; \frac{1}{2}\left(\delta_{\xi_1,\xi_2} - \delta_{\xi_1,-\xi_2}\right), \qquad \xi_1,\xi_2\in 2\pi\Z,
    \\
 \int_0^1\cos(\xi_1 x)\sin(\xi_2 x)dx = &\;0, \qquad \xi_1,\xi_2\in 2\pi\Z,
\end{align*}
to observe that, by the definition of $\Gamma$ and $\ska{\cdot}{\cdot}_{\mathbb A}$,
\begin{align}
S(u,v) &= \ska{\Gamma^0(\tilde u,\tilde v)}{\Gamma^0(\tilde u,\tilde v)}_{\tilde{\mathbb A}}-\ska{\Gamma^0(\tilde u,\tilde u)}{\Gamma^0(\tilde v,\tilde v)}_{\tilde{\mathbb A}}\nonumber\\
&+\frac{1}{4}\int_\S\partial_x(\alpha\cos l_1x+\beta\cos k_1x)A^{-1}\partial_x(\alpha\cos l_1x+\beta\cos k_1x)\,dx\nonumber\\
&-\int_\S\partial_x(\alpha\cos k_1x)A^{-1}\partial_x(\beta\cos l_1x)\,dx\nonumber\\
&-\int_\S\partial_x(\cos k_1x\cos l_1x+\tfrac{1}{2}k_1l_1\sin k_1x\sin l_1x+\tfrac{1}{2}\cos k_2x\cos l_2x)\times\nonumber\\
&\qquad\times A^{-1}\partial_x(\alpha\cos l_1x+\beta\cos k_1x)\,dx\nonumber\\
&+\int_S\partial_x(\cos^2k_1x+\tfrac{1}{2}k_1^2\sin^2k_1x+\tfrac{1}{2}\cos^2k_2x)A^{-1}\partial_x(\beta\cos l_1x)\,dx\nonumber\\
&+\int_S\partial_x(\cos^2l_1x+\tfrac{1}{2}l_1^2\sin^2l_1x+\tfrac{1}{2}\cos^2l_2x)A^{-1}\partial_x(\alpha\cos k_1x)\,dx.\nonumber
\end{align}
By \cite[Prop.\ 5.1]{EKL11} the sum of the first two terms equals the (non-normalized) sectional curvature $S_{\text{2CH}}(\tilde u,\tilde v)$ for the 2CH equation without vorticity and
$$S_{\text{2CH}}\left(\begin{pmatrix}\cos k_1x\\\cos k_2x\end{pmatrix},\begin{pmatrix}\cos l_1x\\\cos l_2x\end{pmatrix}\right)\geq\frac{k_1^2l_1^2}{16}\left[\frac{1}{k_1^2l_1^2} -
\frac{1}{k_1l_1} + \frac{1}{4} - \frac{1}{2 k_1^2 l_1^2} -
\frac{2}{k_1l_1}\right]>0.$$
We may thus write
$$S(u,v)=S_{\text{2CH}}(\tilde u,\tilde v)+\mathcal K_1+\mathcal K_2+\mathcal K_3+\mathcal K_4+\mathcal K_5$$
and conclude that
\begin{align}
\mathcal K_1 &= \frac{1}{4}\left(\frac{1}{2}\alpha^2+\frac{1}{2}\beta^2+\alpha\beta(\delta_{k_1,l_1}+\delta_{k_1,-l_1})\right)\nonumber\\
&-\frac{1}{4}\left(\frac{\alpha^2}{2(1+l_1^2)}+\frac{\beta^2}{2(1+k_1^2)}+\frac{1}{2}\left(\frac{\alpha\beta}{1+l_1^2}
+\frac{\alpha\beta}{1+k_1^2}\right)(\delta_{k_1,l_1}+\delta_{k_1,-l_1})\right),\nonumber\\
\mathcal K_2 &= \frac{1}{2}\alpha\beta\left(-\delta_{k_1,l_1}-\delta_{k_1,-l_1}+\frac{1}{1+l_1^2}(\delta_{k_1,l_1}+\delta_{k_1,-l_1})\right),\nonumber\\
\mathcal K_3 &= -\frac{\alpha}{2}\frac{l_1^2}{1+l_1^2}\bigg[\frac{1}{2}(\delta_{k_1+l_1,l_1}+\delta_{k_1+l_1,-l_1}+\delta_{k_1-l_1,l_1}+\delta_{k_1-l_1,-l_1})\nonumber\\
&+\frac{1}{4}k_1l_1(\delta_{k_1-l_1,l_1}+\delta_{k_1-l_1,-l_1}-\delta_{k_1+l_1,l_1}-\delta_{k_1+l_1,-l_1})\nonumber\\
&+\frac{1}{4}(\delta_{k_2+l_2,l_1}+\delta_{k_2+l_2,-l_1}+\delta_{k_2-l_2,l_1}+\delta_{k_2-l_2,-l_1})\bigg]\nonumber\\
&-\frac{\beta}{2}\frac{k_1^2}{1+k_1^2}\bigg[\frac{1}{2}(\delta_{k_1+l_1,k_1}+\delta_{k_1+l_1,-k_1}+\delta_{k_1-l_1,k_1}+\delta_{k_1-l_1,-k_1})\nonumber\\
&+\frac{1}{4}k_1l_1(\delta_{k_1-l_1,k_1}+\delta_{k_1-l_1,-k_1}-\delta_{k_1+l_1,k_1}-\delta_{k_1+l_1,-k_1})\nonumber\\
&+\frac{1}{4}(\delta_{k_2+l_2,k_1}+\delta_{k_2+l_2,-k_1}+\delta_{k_2-l_2,k_1}+\delta_{k_2-l_2,-k_1})\bigg],\nonumber\\
\mathcal K_4 &=\frac{\beta}{4}\frac{l_1^2}{1+l_1^2}\left((1-k_1^2)(\delta_{2k_1,l_1}+\delta_{2k_1,-l_1})+\delta_{2k_2,l_1}+\delta_{2k_2,-l_1}\right),\nonumber\\
\mathcal K_5 &=\frac{\alpha}{4}\frac{k_1^2}{1+k_1^2}\left((1-l_1^2)(\delta_{2l_1,k_1}+\delta_{2l_1,-k_1})+\delta_{2l_2,k_1}+\delta_{2l_2,-k_1}\right).\nonumber
\end{align}
We let $\alpha>0$ and $\beta=1$ and recall that $k_1,l_1,k_2,l_2\in\{2\pi,4\pi,...\}$ so that at most one Kronecker delta within each pair $\delta_{\xi_1,\xi_2}+\delta_{\xi_1,-\xi_2}$ gives a nonzero contribution. Then a lower estimate for the sectional curvature is given by
\begin{align}
S(u,v)&\geq \frac{1}{8}(\alpha^2+1)-\frac{\alpha^2+1}{8(1+4\pi^2)}-\frac{\alpha}{4(1+4\pi^2)}-\frac{\alpha}{2}-\frac{\alpha+1}{2}
\left(\frac{3}{2}+\frac{1}{4}k_1l_1\right)\nonumber\\
&\quad-\frac{1}{4}\frac{k_1^2l_1^2}{1+l_1^2}-\frac{\alpha}{4}\frac{k_1^2l_1^2}{1+k_1^2}\nonumber\\
&\geq\frac{1}{10}(\alpha^2-5\alpha M^2-5M^2),\nonumber
\end{align}
where $M=\max\{k_1,l_1\}$. The right hand side of the above inequality is positive for
$$\alpha\geq 6M^2>\frac{5}{2}M^2+\sqrt{\frac{25}{4}M^4+5M^2}.$$
As
\begin{align}
\ska{u}{u}_{\mathbb A}&=1+\tfrac{1}{2}(k_1^2+\alpha^2),\nonumber\\
\ska{v}{v}_{\mathbb A}&=1+\tfrac{1}{2}(l_1^2+1)\text{ and}\nonumber\\
\ska{u}{v}_{\mathbb A}&=\tfrac{\alpha}{2},\nonumber
\end{align}
it is clear that
$$K(u,v)\geq\frac{\frac{1}{10}(\alpha^2-5\alpha M^2-5M^2)}{(1+\frac{1}{2}(k_1^2+\alpha^2))(1+\frac{1}{2}(l_1^2+1))-\frac{\alpha^2}{4}}\to
\frac{\tfrac{1}{10}}{\tfrac{1}{2}(1+\frac{1}{2}(l_1^2+1))-\frac{1}{4}}>0$$
as $\alpha\to\infty$. Thus the proof of Theorem~\ref{thm_curvature} is completed.\hfill$\square$
\end{document}